\documentclass[12pt]{article}

\title{\bf On Amicable Numbers With Different Parity}

\author{Germano D'Abramo\\
{\small Istituto di Astrofisica Spaziale e Fisica Cosmica, Roma.}\\
{\small E--mail: {\tt Germano.Dabramo@rm.iasf.cnr.it}}}

\date{}

\begin{document}

%

\maketitle

\abstract{In this paper we provide a straightforward proof that if a pair 
of amicable numbers with different parity exists (one number odd and the 
other one even), then the odd amicable number must be a perfect square, 
while the even amicable number has to be equal to the product of a power 
of $2$ and an odd perfect square.}

\section{Introductory remarks}

An amicable pair $(M,N)$ consists of two integers $M$,~$N$ for which the 
sum of proper divisors (the divisors excluding the number itself) of one 
number equals the other.

Amicable numbers were known to the Pythagoreans and many outstanding 
mathematicians of the past (Pierre de Fermat, Ren\'e Descartes and 
Leonhard Euler to quote a few) were interested in them and discovered new 
couples~\cite{bi1}.

Some algorithms are known to find some kind of amicable couples (that of 
the Arab mathematician Thabit~ibn~Qurra, for instance, subsequently 
extended by Euler's~algorithm), but no one has yet found a single method 
that generates all possible amicable pairs. Nowadays mathematicians know 
of more than 1~million amicable pairs~\cite{bi1,bi2}.

It is not known whether there exist infinitely many pairs of amicable 
numbers. Another unsolved puzzle concerns the parity of amicable pairs. In 
every known pair, both numbers are even or both are odd. However, no one 
has yet proved that no pair exists in which one number is odd and the 
other is even~\cite{bi1}.

Suppose you want to prove that a pair of amicable numbers with different 
parity cannot exist. A possible approach (a trivial one, in fact) is to 
prove that the existence of the even number of the pair is incompatible 
with both the property of being odd of the other number belonging to the 
pair and the mathematical property of being `amicable' of the pair.

For the sake of the argument, in this paper we try to proceed along this 
direction and see where it leads, namely which interesting properties can 
be deduced about possible existing pairs of amicable numbers with 
different parity.

\section{Properties of the odd amicable number}

Let us assume that $N$ is an odd number belonging to a pair of amicable 
numbers. In what follows we study the parity of the second number $M$ of 
the pair, depending on a few simple number-theoretical properties of $N$.

In order for $N$ to be an odd number it is necessary that its prime-power 
factorization is of the kind

\begin{equation}
N=p^{k_{1}}_1\cdot p^{k_{2}}_2\cdot p^{k_{3}}_3\cdot ...\cdot p^{k_{n}}_n,
\label{eq1}
\end{equation}
where the $p_i$'s are all {\em odd} prime numbers (i.e.~primes different 
from $2$) and the $k_i$'s are integer numbers.

From equation~(\ref{eq1}) it turns out to be that all the proper divisors 
of $N$ are odd, since they are either odd prime numbers or products of 
powers of odd prime numbers. So, in order to study the parity of the 
amicable number associated to $N$, it is sufficient to study the parity of 
the sum of all the proper divisors of $N$. As a matter of fact, such sum 
will be even only if $N$ has an even number of proper divisors, otherwise 
the sum will be odd.

In the following three subsections we show how the parity of the total 
number of proper divisors of $N$ behaves according to how $N$ is 
factorized in the form of equation~(\ref{eq1}).

\subsection*{Case 1}

Let us start with $N=p^{k_{1}}_1$. The set of its proper divisors is then

\begin{equation}
1,\,p_1,\,p^{2}_1,\, ...\, ,\, p^{k_{1}-1}_1,
\label{eq2}
\end{equation}
and consequently the second amicable number is the sum of $k_{1}$ odd 
numbers. It follows that the second amicable number is even only if 
$k_{1}$ is an even number, namely only if $N=p^{k_{1}}_1$ is a perfect 
square.

\subsection*{Case 2}

Let us have now $N=p^{k_{1}}_1\cdot p^{k_{2}}_2$. The sum of all its 
proper divisors, that we call $S_2$, is

\begin{displaymath}
S_2=(p_1+p^{2}_1+\cdots +p^{k_{1}}_1)+(p_2+p^{2}_2+\cdots +p^{k_{2}}_2)+
p^{k_{2}}_2\cdot(p_1+p^{2}_1+\cdots +p^{k_{1}-1}_1)+
\end{displaymath}
\begin{equation}
+(p_1+p^{2}_1+\cdots +p^{k_{1}}_1)\cdot(p_2+p^{2}_2+\cdots +p^{k_{2}-1}_2)+1.
\label{eq3}
\end{equation}

If $k_1-1=0$ or $k_2-1=0$, the terms $p^{k_{2}}_2\cdot(p_1+p^{2}_1+\cdots 
+p^{k_{1}-1}_1)$ or $(p_1+p^{2}_1+\cdots 
+p^{k_{1}}_1)\cdot(p_2+p^{2}_2+\cdots +p^{k_{2}-1}_2)$ must be canceled in 
equation~(\ref{eq3}).

Now, if we develop all the products and count the addends in the second 
member of equation~(\ref{eq3}), it is easy to verify that the total number 
of proper divisors (i.e.~the cardinality of the set of all the addends) is 
$n_2=k_1+k_2+k_1\cdot k_2$. This relation holds also when $k_1=1$ and 
$k_2\geq 1$ (likewise for $k_2=1$ and $k_1\geq 1$).

Note that if $k_1$ and $k_2$ are both even numbers, then $n_2$ is even and 
the sum $S_2$, namely the second amicable number, is an even number too.

If, otherwise, $k_1$ or $k_2$ is an odd number, then $n_2$ is odd and the 
second amicable number is odd too.

So, even in this case we could have two amicable number with different 
parity only if $N=p^{k_{1}}_1\cdot p^{k_{2}}_2$ with both $k_1$ and $k_2$ 
even, namely only if it is a perfect square.

\subsection*{Case 2: Extension}

Let $N=p^{k_{1}}_1\cdot p^{k_{2}}_2\cdot p^{k_{3}}_3$. We can use the 
results of Case~2 to derive the sum of all the proper divisors of $N$ and 
study its parity.

The sum of all the proper divisors of $N=p^{k_{1}}_1\cdot p^{k_{2}}_2\cdot 
p^{k_{3}}_3$ is

\begin{equation}
S_2\cdot (1+p_3+p^{2}_3+\cdots +p^{k_{3}}_3)+
p^{k_{1}}_1\cdot p^{k_{2}}_2\cdot(1+p_3+p^{2}_3+\cdots +p^{k_{3}-1}_3),
\label{eq4}
\end{equation}
where $S_2$ is the sum shown in equation~(\ref{eq3}). If $k_3-1=0$ we must 
put $p^{k_{1}}_1\cdot p^{k_{2}}_2$ instead of $p^{k_{1}}_1\cdot 
p^{k_{2}}_2\cdot(1+p_3+p^{2}_3+\cdots +p^{k_{3}-1}_3)$.

Also in this case, once carried out all the products, it is easy to see 
that the total number of proper divisors of $N$ is $n_3=k_3\cdot 
(n_2+1)+n_2$, where $n_2$ is the total number of proper divisors derived 
in Case~2 (which, we remind, is $k_1+k_2+k_1\cdot k_2$). Note that 
relation $k_3\cdot (n_2+1)+n_2$ holds also for $k_3=1$.

Let us now study the parity. If $n_2$ is odd, then $n_3$ is odd too, 
independently of the parity of $k_3$. According to the results of Case~2, 
in order for $n_2$ to be odd it is necessary that $k_1$ or $k_2$ is an odd 
number. Therefore, it is necessary that $p^{k_{1}}_1\cdot p^{k_{2}}_2$ is 
{\em not} a perfect square. If the previous condition is verified, then 
the second amicable number will have the same parity of $N$ (i.e. it will 
be odd), independently of the parity of $k_3$.

At this point, the above procedure can be iterated with further prime 
factors $p_i$ in the factorization of $N$.

\subsection*{Generalization}

Given the general prime factorization $N=p^{k_{1}}_1\cdot p^{k_{2}}_2\cdot 
p^{k_{3}}_3\cdot ...\cdot p^{k_{n}}_n$, it is easy to verify that the 
number $n_n$ of all the proper divisors of $N$ is obtained iteratively as 
$n_i=k_i\cdot (n_{i-1}+1)+n_{i-1}$, for~$i=2,3,4,...,n$. We remind that 
$n_1=k_1$, as shown in Case~1.

Since the choice of the first prime factors to study in the same order as 
in Case~1 and Case~2 is free, we are allowed to conclude that only if {\em 
all} the $k_i$'s are even then an even amicable number can exist 
associated to the odd number $N$.

As a matter of fact, it is sufficient to have a single prime factor raised 
to an odd power (for example, $p^{k_d}_d$) to apply the arguments showed 
in Case~2, namely taking that factor, multiplying it by any other factor 
(even one raised to an even power), for example $p^{k_d}_d\cdot 
p^{k_j}_j$, and starting the analysis done in Case~2. In this case $n_2$ 
results to be odd and, according to equation $n_i=k_i\cdot 
(n_{i-1}+1)+n_{i-1}$, every $n_i$ will be odd for all the $i$'s. As a 
trivial consequence, $n_n$ will be odd too, implying that the second 
amicable number $M$ will be odd. Therefore, the only possibility for $M$ 
to be even is that all the $k_i$'s have to be even, namely that $N$ must 
be a perfect square.

\section{Properties of the even amicable number}

Let us now consider some properties of the even number $M$ belonging to a 
pair of amicable number with different parity.

The number $M$ can be uniquely factorized in the following way

\begin{equation}
M=2^{k_0}\cdot p^{k_{1}}_1\cdot p^{k_{2}}_2\cdot p^{k_{3}}_3\cdot ...
\cdot p^{k_{n}}_n,
\label{eq5}
\end{equation}
where, again, the $p_i$'s are all {\em odd} prime numbers (i.e.~primes 
different from $2$) and the $k_i$'s are integer numbers (with $k_0\neq 
0$).

If $M$ is simply equal to $2^{k_0}$, namely if $\forall i\neq 0\,\,\, 
k_i=0$, then the sum of all its proper divisors (i.e.~$N$) is

\begin{equation}
N=1+2+2^2+\cdots +2^{k_{0}-1}=2^{k_0}-1.
\label{eq6}
\end{equation}

This number is obviously odd and the relation~(\ref{eq6}) does not violate 
the starting assumption that $N$ is odd. But, as shown in Section~2, $N$ 
has to be also a perfect square, so let us see if and when $2^{k_0}-1$ is 
a perfect square.

In order for $2^{k_0}-1$ to be an odd perfect square we must have that

\begin{equation}
2^{k_0}=(2j+1)^2+1,
\label{eq7}
\end{equation}
with $j\in \{0,1,2,3,...\}$. Therefore, assume that $\forall k_0$ it is 
always possible to find a $j$ such that equation~(\ref{eq7}) holds. Then 
we have

\begin{equation}
2^{k_0}=4j^2+4j+2.
\label{eq8}
\end{equation}

Dividing both members of equation~(\ref{eq8}) by~$2$, we have

\begin{equation}
2^{k_{0}-1}=2(j^2+j)+1,
\label{eq9}
\end{equation}
namely we obtain the paradoxical situation that the first member 
$2^{k_{0}-1}$ is even, while the second one is clearly odd. The only 
possibility for equation~(\ref{eq9}) to be valid is with $k_0$~equal 
to~$1$, namely with $M=2$ and $N=1$, but this pair is {\bf not} a pair of 
amicable numbers.

Summing up briefly, in order for $M$ to be an even number belonging to a 
pair of amicable numbers with different parity it is necessary that it is 
of the form~(\ref{eq5}) with {\em at least} one {\em odd} prime factor.

Now we derive the sum of all the proper divisors of $M$. If with $S_d$ we 
call the sum of all the proper divisors of the odd portion of 
factorization~(\ref{eq5}), namely of $p^{k_{1}}_1\cdot p^{k_{2}}_2\cdot 
p^{k_{3}}_3\cdot ...\cdot p^{k_{n}}_n$, then it is easy to check that the 
sum of all the proper divisors of $M$, namely the odd amicable number $N$, 
is

\begin{equation}
N=(1+2+2^2+\cdots +2^{k_0})\cdot S_d + 
(p^{k_{1}}_1\cdot p^{k_{2}}_2\cdot p^{k_{3}}_3\cdot ...
\cdot p^{k_{n}}_n)\cdot (1+2+2^2+\cdots +2^{{k_{0}-1}}),
\label{eq10}
\end{equation}
and, rewriting all the sums of powers of $2$ in the following compact way, 
we have

\begin{equation}
N=(2^{k_{0}+1}-1)\cdot S_d + (p^{k_{1}}_1\cdot p^{k_{2}}_2\cdot 
p^{k_{3}}_3\cdot ...\cdot p^{k_{n}}_n)\cdot (2^{k_{0}}-1).
\label{eq11}
\end{equation}

It is easy to note that the addend $(p^{k_{1}}_1\cdot p^{k_{2}}_2\cdot 
p^{k_{3}}_3\cdot ...\cdot p^{k_{n}}_n)\cdot (2^{k_{0}}-1)$ of 
equation~(\ref{eq11}) is an odd number (since it is a product of odd 
numbers), as well as the term $(2^{k_{0}+1}-1)$ in the first addend.

Therefore, in order for $N$ to be odd, as given by assumption, it is 
necessary that $S_d$, namely the sum of all the proper divisors of the odd 
portion of factorization~(\ref{eq5}), has to be even. But, according to 
what has been proved in Section~2, this is possible only if the odd 
portion $p^{k_{1}}_1\cdot p^{k_{2}}_2\cdot p^{k_{3}}_3\cdot ...\cdot 
p^{k_{n}}_n$ of the factorization of $M$ is a perfect square.

At this point we have proved our theorem.

\section{Conclusion}

In the present note we have proved that if a pair of amicable numbers with 
different parity exists, then the odd number must be a perfect square, 
while the even number has to be equal to the product of a power of $2$ and 
an odd perfect square.

This theorem might be useful in the implementation of an algorithm for 
numerical search of possible existing pairs of amicable numbers with 
different parity. Hopefully, it might turn out to be useful also within a 
future, wider theorem which proves the (non)existence of such pairs.

\end{document}